\documentclass{article}

\input{tcilatex}

\begin{document}

\title{The Use of Geometric Quantities in the Tensor Description of a Euclidean Space}
\author{Pavel Grinfeld}
\maketitle
\begin{abstract}
We present a tensor description of Euclidean spaces that emphasizes the use of geometric vectors. We demonstrate the
effectiveness of the approach by proving of a number of integral identities with vector integrands. Visit grinfeld.org for the 
latest version.
\end{abstract}%

\section{Introduction}

Since the invention of coordinate systems in the middle of the 17th century,
the subject of Geometry has followed the steady path of algebraization. This
should not surprise us: the very idea of a coordinate system is to replace
geometric objects with the coordinates of the constituent points, thus
opening the problem up to algebraic -- and, since the invention of Calculus, 
\textit{analytical} -- methods. This gives the method of coordinates a
distinct advantage over its geometric peers: while geometric arguments
typically require a unique insight into the problem, and therefore a certain
degree of ingenuity, algebraic and analytical methods tend towards
universality and robustness. The advent of computing has further cemented
this advantage.

On the other hand, the use of coordinates comes with great costs. Chief
among them is the loss of geometric insight. This pitfall is exemplified by
the results of Leonhard Euler and Louis Lagrange in their foundational works
on the Calculus of Variations. In $1744$, in his search for a minimal shape
of revolution, Euler introduced what we would now call a cylindrical
coordinate system and described the profile of the minimal surface by an
unknown function $r\left( z\right) $, as illustrated in the following figure.%
\begin{equation}
\FRAME{itbpF}{1.4166in}{1.5359in}{0in}{}{}{catenoidwithmarkings.png}{\special%
{language "Scientific Word";type "GRAPHIC";maintain-aspect-ratio
TRUE;display "USEDEF";valid_file "F";width 1.4166in;height 1.5359in;depth
0in;original-width 13.8889in;original-height 14.6248in;cropleft "0";croptop
"1";cropright "1";cropbottom "0";filename
'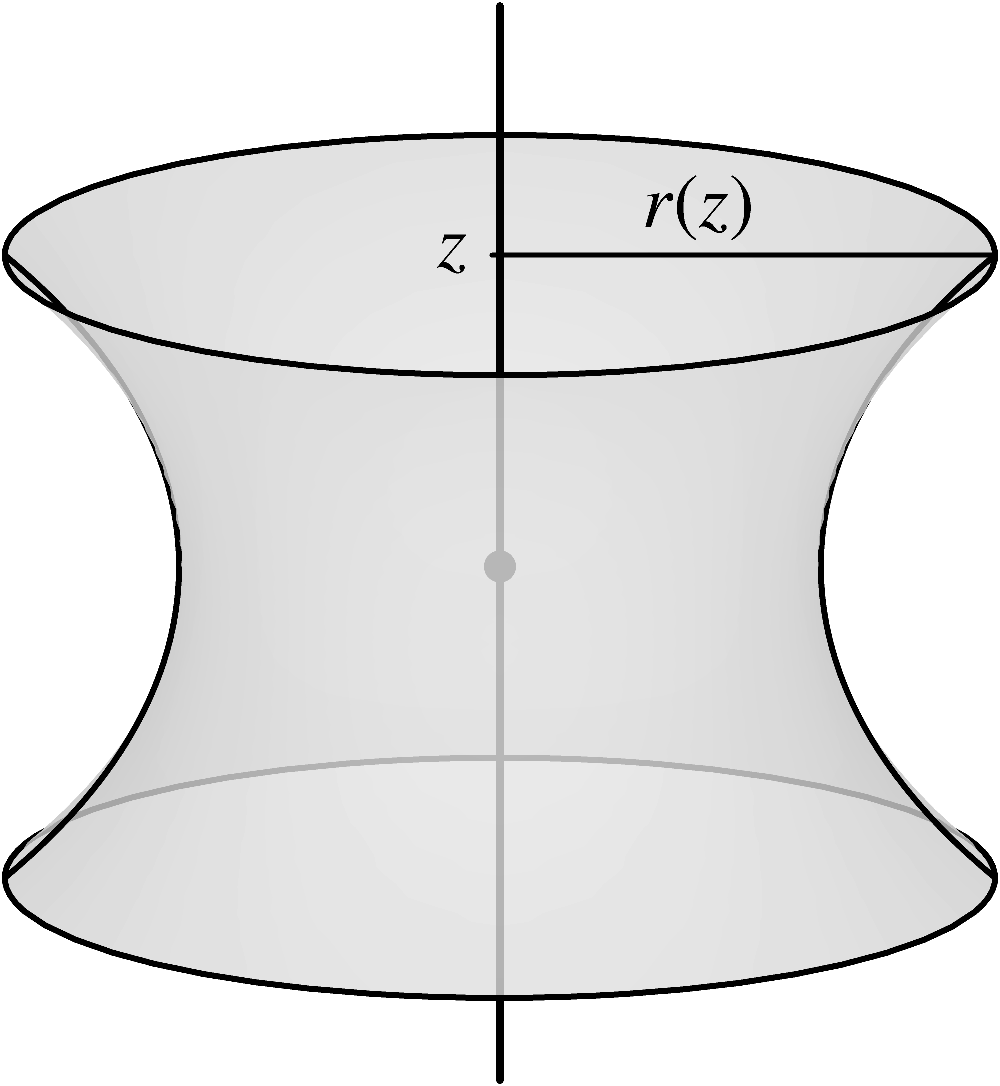';file-properties "XNPEU";}}
\end{equation}%
By making arguments based on manipulating \textit{geometric} elements, Euler
demonstrated that $r\left( z\right) $ must satisfy the equation%
\begin{equation}
r^{\prime \prime }\left( z\right) r\left( z\right) -r^{\prime }\left(
z\right) ^{2}-1=0.  \label{Catenoid eqn}
\end{equation}%
Euler then promptly solved this equation to reveal that 
\begin{equation}
r\left( z\right) =a\cosh \frac{z}{a}.
\end{equation}%
In other words, the profile is a \textit{catenary}, where the constant $a$
represents the closest distance between the surface and the axis of
revolution.

In $1755$, the nineteen-year-old Lagrange took an even more unapologetically
coordinate approach to the problem of minimal surfaces and represented the
unknown surface by the graph of a function%
\begin{equation}
z=F\left( x,y\right)
\end{equation}%
in Cartesian coordinates $x,y,z$.%
\begin{equation}
\FRAME{itbpF}{2.8055in}{1.4823in}{0in}{}{}{lagrangesurface.png}{\special%
{language "Scientific Word";type "GRAPHIC";maintain-aspect-ratio
TRUE;display "USEDEF";valid_file "F";width 2.8055in;height 1.4823in;depth
0in;original-width 13.8889in;original-height 7.2774in;cropleft "0";croptop
"1";cropright "1";cropbottom "0";filename
'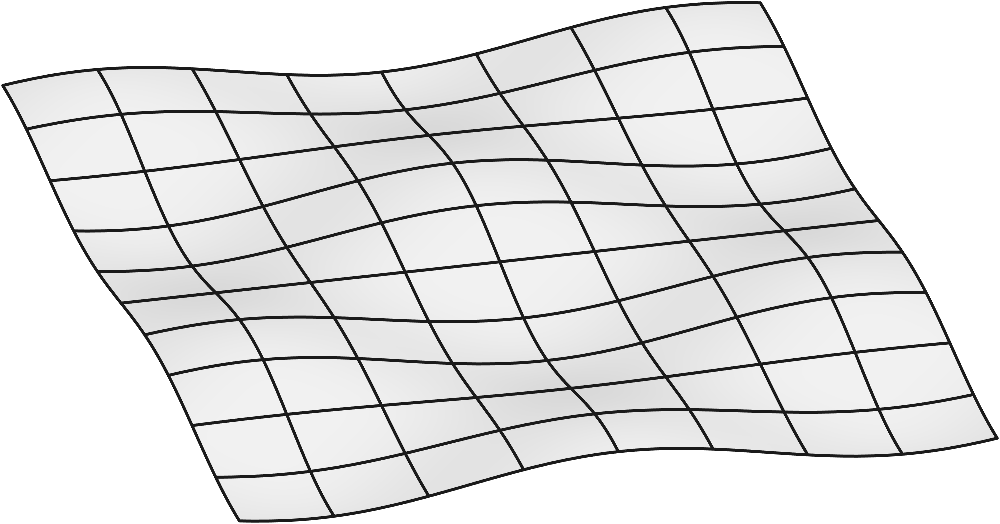';file-properties "XNPEU";}}
\end{equation}%
Reasoning \textit{analytically}, Lagrange demonstrated that $z\left(
x,y\right) $ must satisfy the partial differential equation%
\begin{equation}
\frac{\partial ^{2}F}{\partial x^{2}}+\frac{\partial ^{2}F}{\partial y^{2}}+%
\frac{\partial ^{2}F}{\partial x^{2}}\left( \frac{\partial F}{\partial y}%
\right) ^{2}+\frac{\partial ^{2}F}{\partial y^{2}}\left( \frac{\partial F}{%
\partial x}\right) ^{2}-2\frac{\partial F}{\partial x}\frac{\partial F}{%
\partial y}\frac{\partial ^{2}F}{\partial x\partial y}=0.
\end{equation}

These works of Euler and Lagrange have been rightfully revered by later
mathematicians for the seminal nature of the techniques developed in them
and for the flight of creativity that their development required. However,
to us, as we contemplate the advantages and disadvantages of coordinate
systems, these works offer an additional insight that speaks to the
potential loss of geometric insight that comes with the use of coordinates.
It appears that neither Euler nor Lagrange found the geometric
interpretation for their equations. It was only in $1774$ that the French
mathematician Jean Baptiste Meusnier discovered that \textit{a minimal
surface is characterized by zero mean curvature}. Critically, the loss of
geometric insight is not some aesthetic lament: geometric insight, after
all, serves as a guide for organizing analytical expressions into meaningful
combinations. As a result, with the loss of geometric insight, we lose
control over our calculations. The algebraic complexity in the analysis
grows rapidly with every step until we are forced to retreat in the face of
computational difficulties.

When it comes to combatting the loss of geometric insight, two distinct
approaches have been developed:\ the tensor calculus approach and the dyadic
approach associated with modern differential geometry. In tensor calculus, a
geometric problem is analyzed by introducing a coordinate system and
replacing all geometric objects by their components. Tensor calculus
preserves the geometric insight by developing a framework of \textit{%
invariance} that dictates precise rules for combing analytical expressions
into geometrically meaningful combinations. By contrast, the dyadic approach
eschews components altogether and operates only in terms of invariant
objects and operators. Naturally, both approaches have their uses and
misuses, their advantages and disadvantages, and their adherents and
detractors. And, as always, the truth is that elements of both approaches
are essential and the two schools of thoughts are complementary rather than
in conflict.

The goal of this paper is to describe a \textit{vector tensor calculus},
i.e. a particular style of tensor treatment of a Euclidean space that
combines elements of both the tensor and dyadic approaches by emphasizing
the use of geometric vectors. Some elements of this approach can be found in
V.F. Kagan's \textit{Foundations of the Theory of Surfaces in Tensor Terms 
\cite{KaganTheoryOfSurfaces}}. However, Kagan typically uses geometric
vector quantities only at the outset of any particular discussion only to
abandon them in favor of working with scalar quantities. This is
understandable -- scalar quantities remain meaningful in the generalization
to Riemannian spaces and are therefore more robust than geometric vector
quantities in this sense. We take the analysis of vector quantities much
further and discover their tremendous utility in simplifying the description
of Euclidean spaces as well as revealing new and insightful relationships.

Another noteworthy aspect of Kagan's textbook that deserves to be mentioned
is its penchant for omitting technical details in favor of transparency when
communicating essential ideas. Interesting, and relevant to the goals of
this paper, is the fact that Kagan's approach did not sit right with many of
his contemporaries. In an otherwise positive review \cite{AlexandrovOnKagan}%
, A.D. Alexandrov criticized Kagan's lack of formalism:

\textit{...Other shortcomings of this book have to do with the prevalence of
the tensor framework. They manifest themselves in insufficient attention to
the precise definitions of concepts and to the specification of assumptions
required for correctness of theorems.}

\textit{It is not my goal to criticize the work of V.F. Kagan. Such
deficiencies are characteristic of an entire direction in differential
geometry and can be found in the majority of books devoted to this field.
They have become a matter of style that I\ find anachronistic, as our
present notion of rigor is different from that of, say, the middle of the
nineteenth century.}

With Alexandrov's remarks duly noted and in the spirit of Kagan's classic,
we too will favor clarity over rigor. We will generally assume that all
surfaces are infinitely smooth. The key takeaway will be that an emphasis on
geometric vectors provides greater insight into the structure of Euclidean
spaces, offers more elegant demonstrations of known results, and opens doors
to new results. For the more standard approach to the tensor description of
Euclidean spaces, see the classical textbooks \cite{LeviCivitaTensors}, \cite%
{McConnellTensors}, \cite{SyngeTensorCalculus}.

\section{Summary of demonstrated identities}

As an illustration of the effectiveness of the proposed approach, we will
demonstrate a family of integral relationships for a smooth closed
hypersurface $S$ with unit normal $\mathbf{N}$, mean curvature $B_{\alpha
}^{\alpha }$, and Gaussian curvature $K$ in an $n$-dimensional Euclidean
space. Naturally, the results involving the Gaussian curvature are limited
to $n=3$. For higher-dimensional spaces, the scalar curvature $R$ takes the
place of the Gaussian curvature $K$.

First, we will prove the well-known fact that the surface integral of the
unit normal vanishes, i.e.%
\begin{equation}
\int_{S}\mathbf{N}dS=\mathbf{0.}  \label{IN = 0}
\end{equation}%
Similarly, we will show that the surface integral of the combination $%
\mathbf{N}B_{\alpha }^{\alpha }$, known as the curvature normal, also
vanishes, i.e.%
\begin{equation}
\int_{S}\mathbf{N}B_{\alpha }^{\alpha }dS=\mathbf{0.}  \label{INH = 0}
\end{equation}%
Finally, we will demonstrate that the surface integral of the combination $%
\mathbf{N}K$ vanishes as well, i.e. 
\begin{equation}
\int_{S}\mathbf{N}KdS=\mathbf{0.}  \label{INK = 0}
\end{equation}

Let $\mathbf{R}$ be the position vector emanating from an arbitrary origin $%
O $. For any vector quantity $\mathbf{U}$ whose surface integral vanishes,
i.e.%
\begin{equation}
\int_{S}\mathbf{U~}dS=\mathbf{0},
\end{equation}%
it is natural to inquire as to the value of the integral%
\begin{equation}
\int_{S}\mathbf{R}\cdot \mathbf{U~}dS
\end{equation}%
since it is independent of the arbitrary origin $O$. After all, for 
\begin{equation}
\mathbf{R}^{\prime }=\mathbf{R}+\mathbf{d},
\end{equation}%
we have 
\begin{eqnarray}
\int_{S}\mathbf{R}^{\prime }\cdot \mathbf{U}dS &=&\int_{S}\left( \mathbf{R}+%
\mathbf{d}\right) \cdot \mathbf{U}dS \\
&=&\int_{S}\mathbf{R}\cdot \mathbf{U}dS+\mathbf{d}\cdot \int_{S}\mathbf{U}dS
\\
&=&\int_{S}\mathbf{R}\cdot \mathbf{U}dS.
\end{eqnarray}%
Independence from $O$ suggests that the integral%
\begin{equation}
\int_{S}\mathbf{R}\cdot \mathbf{U}dS
\end{equation}%
represents a geometric characteristic of the surface $S$. Indeed, for each
of the vector fields $\mathbf{N}$, $\mathbf{N}B_{\alpha }^{\alpha }$, and $%
\mathbf{N}K$, the surface integral of the dot product with the position
vector $\mathbf{R}$ yields an interested geometric quantity. Namely,%
\begin{eqnarray}
\int_{S}\mathbf{R}\cdot \mathbf{N}dS\ \ \ ~ &=&nV  \label{IR.N} \\
\int_{S}\mathbf{R}\cdot \mathbf{N}B_{\alpha }^{\alpha }dS &=&-\left(
n-1\right) A  \label{IR.NH} \\
\int_{S}\mathbf{R}\cdot \mathbf{N}KdS\ \ &=&-\frac{1}{2}\int_{S}B_{\alpha
}^{\alpha }dS,  \label{IR.NK}
\end{eqnarray}%
where $V$ is the volume of enclosed domain and $A$ is the surface area of $S$%
.

Note that the same logic applies to the integral%
\begin{equation}
\int_{S}\mathbf{R}\times \mathbf{U}dS
\end{equation}%
and we would find that%
\begin{equation}
\int_{S}\mathbf{R}\times \mathbf{N}dS=\int_{S}\mathbf{R}\times \mathbf{N}%
B_{\alpha }^{\alpha }dS=\int_{S}\mathbf{R}\times \mathbf{N}KdS=\mathbf{0}.
\end{equation}

\section{A tensor description of a Euclidean space}

In this Section, we will describe the fundamental tensor objects in a
Euclidean space. As we have already mentioned above, the distinguishing
characteristic of our description is its emphasis on geometric quantities.
Since it is not possible to present a full account in the limited space, we
will only give the definitions of the key objects, state their fundamental
properties, and list the essential identities relating those objects. A
detailed description of the presented approach can be found in \cite%
{GrinfeldTC}.

Refer the ambient Euclidean space to arbitrary curvilinear coordinates $%
Z^{1},Z^{2},Z^{3}$ or, collectively, $Z^{i}$, and treat the position vector $%
\mathbf{R}$ as a function of $Z^{i}$, i.e.%
\begin{equation}
\mathbf{R}=\mathbf{R}\left( Z\right) ,
\end{equation}%
where the shorthand symbol $\mathbf{R}$ represents the function $\mathbf{R}%
\left( Z^{1},Z^{2},Z^{3}\right) $. Then \textit{covariant basis} $\mathbf{Z}%
_{i}$, the \textit{contravariant basis} $\mathbf{Z}^{i}$, the \textit{%
covariant metric tensor} $Z_{ij}$, and the \textit{contravariant metric
tensor} $Z^{ij}$ are given by the identities%
\begin{eqnarray}
\mathbf{Z}_{i} &=&\frac{\mathbf{\partial R}\left( Z\right) }{\partial Z^{i}}
\label{Zi} \\
\mathbf{Z}^{i}\cdot \mathbf{Z}_{j} &=&\delta _{j}^{i}  \label{ZI.Zj} \\
Z_{ij} &=&\mathbf{Z}_{i}\cdot \mathbf{Z}_{j}  \label{Zij} \\
Z^{ij}Z_{jk} &=&\delta _{k}^{i},  \label{ZIJ}
\end{eqnarray}%
where $\delta _{k}^{i}$ is the familiar \textit{Kronecker delta symbol}.

Suppose that the surface $S$ is referred to the surface coordinates $%
S^{1},S^{2}$ or, collectively $S^{\alpha }$, and treat the surface
restriction of the position vector $\mathbf{R}$ as a function of $S^{\alpha
} $, i.e.%
\begin{equation}
\mathbf{R}=\mathbf{R}\left( S\right) .
\end{equation}%
Then \textit{covariant basis} $\mathbf{S}_{\alpha }$, the \textit{%
contravariant basis} $\mathbf{S}^{\alpha }$, the \textit{covariant metric
tensor} $S_{\alpha \beta }$, and the \textit{contravariant metric tensor} $%
S^{\alpha \beta }$ are given by the identities%
\begin{eqnarray}
\mathbf{S}_{\alpha } &=&\frac{\mathbf{\partial R}\left( S\right) }{\partial
S^{\alpha }}  \label{Sa} \\
\mathbf{S}^{\alpha }\cdot \mathbf{S}_{\beta } &=&\delta _{\beta }^{\alpha }
\label{SA.Sb} \\
S_{\alpha \beta } &=&\mathbf{S}_{\alpha }\cdot \mathbf{S}_{\beta }
\label{Sab} \\
S^{\alpha \beta }S_{\beta \gamma } &=&\delta _{\gamma }^{\alpha }.
\label{SAB}
\end{eqnarray}

The components $U^{i}$ of a vector $\mathbf{U}$ are given by the dot product
of $\mathbf{U}$ with the contravariant basis $\mathbf{Z}^{i}$, i.e.%
\begin{equation}
U^{i}=\mathbf{Z}^{i}\cdot \mathbf{U.}
\end{equation}%
Similarly, the surface components $U^{\alpha }$ of a vector $\mathbf{U}$ is
the plane tangential to surface $S$ are given by the dot product of $\mathbf{%
U}$ with the surface contravariant basis $\mathbf{S}^{\alpha }$, i.e.%
\begin{equation}
U^{\alpha }=\mathbf{S}^{\alpha }\cdot \mathbf{U.}
\end{equation}%
The shift tensor $Z_{\alpha }^{i}$ represents the ambient coordinates of the
surface covariant basis $\mathbf{S}_{\alpha }$, i.e.%
\begin{equation}
Z_{\alpha }^{i}=\mathbf{Z}^{i}\cdot \mathbf{S}_{\alpha }.
\end{equation}%
The shift tensor is a critical object in the traditional approach to tensor
Calculus, but will not figure in our analysis.

For a two-dimensional hypersurface, the unit normal $\mathbf{N}$ is given by
the identity%
\begin{equation}
\mathbf{N}=\frac{1}{2}\varepsilon ^{\alpha \beta }\mathbf{S}_{\alpha }\times 
\mathbf{S}_{\beta },  \label{N}
\end{equation}%
where $\varepsilon ^{\alpha \beta }$ is the surface \textit{Levi-Civita
symbol}. It is a straightforward matter to generalize this identity to
higher dimensions.

The \textit{covariant derivative} $\nabla _{k}$ is a differential operator
that preserves the tensor property of its inputs. It satisfies product rule,
the sum rule, and the metrinilic property with respect to all the ambient
metrics, i.e. 
\begin{eqnarray}
\nabla _{k}Z_{ij},\nabla _{k}\delta _{j}^{i},\nabla _{k}Z^{ij},\nabla
_{k}\varepsilon _{rst} &=&0  \label{MN1} \\
\nabla _{k}\mathbf{Z}_{i},\nabla _{k}\mathbf{Z}^{i} &=&\mathbf{0}
\label{MN2}
\end{eqnarray}%
Furthermore, the covariant derivative $\nabla _{k}$ coincides with the
partial derivative $\partial /\partial Z^{k}$ in affine coordinates as well
as for tensors of order zero in arbitrary coordinates. In particular, the
covariant basis $\mathbf{Z}_{i}$ can be expressed in terms of the covariant
derivative, i.e.%
\begin{equation}
\mathbf{Z}_{i}=\nabla _{i}\mathbf{R.}  \label{Zi = ,iR}
\end{equation}

The \textit{surface covariant derivative} $\nabla _{\gamma }$ is a
differential operator that applies to objects defined on the surface. It is,
too, distinguished by the property that it preserves the tensor property of
its inputs. It satisfies the sum rule, the product rule, and the metrinilic
property with respect to the metric tensors and the Levi-Civita symbols, i.e.%
\begin{equation}
\nabla _{\gamma }S_{\alpha \beta },\nabla _{\gamma }\delta _{\beta }^{\alpha
},\nabla _{\gamma }S^{\alpha \beta }=0,
\end{equation}%
but not with respect to the surface bases $\mathbf{S}_{\alpha }$ and $%
\mathbf{S}^{\alpha }$. It coincides with the partial derivative $\partial
/\partial S^{\gamma }$ in affine coordinates, provided that the surface
admits such coordinates. It also coincides with the partial derivative for
tensors of order zero in arbitrary coordinates. In particular,%
\begin{equation}
\mathbf{S}_{\alpha }=\nabla _{\alpha }\mathbf{R.}  \label{Sa = ,aR}
\end{equation}

The ambient version of the \textit{divergence theorem} states that the
volume integral of an invariant quantity $\nabla _{i}T^{i}$ equals the
surface integral of the invariant quantity $N_{i}T^{i}$, i.e.%
\begin{equation}
\int_{\Omega }\nabla _{i}T^{i}dZ=\int_{S}N_{i}T^{i}dS,  \label{DT}
\end{equation}%
where $\Omega $ is the domain enclosed by the surface $S$, and $N_{i}$ are
the components of the \textit{external} normal $\mathbf{N}$. The theorem is
valid for a tensor $T^{i}$ with either scalar or vector elements.

Similarly, according to the \textit{surface} divergence theorem for a patch $%
S$ with a boundary $L$, the surface integral of the invariant quantity $%
\nabla _{\alpha }T^{\alpha }$ equals the boundary integral of the invariant $%
n_{\alpha }T^{\alpha }$, i.e.%
\begin{equation}
\int_{S}\nabla _{\alpha }T^{\alpha }dS=\int_{L}n_{\alpha }T^{\alpha }dS,
\label{SDT}
\end{equation}%
where $n_{\alpha }$ are the surface components of the exterior unit vector $%
\mathbf{n}$ that is normal to the boundary $L$ and lies in the plane tangent
to the surface. Once again, the statement is valid for a tensor $T^{\alpha }$
with either scalar or vector elements.

The concept of curvature arises in the analysis of the vectors $\nabla
_{\alpha }\mathbf{S}_{\beta }$. While $\nabla _{\alpha }\mathbf{S}_{\beta }$
do not vanish, they are shown to be orthogonal to the surface, i.e.%
\begin{equation}
\nabla _{\alpha }\mathbf{S}_{\beta }=\mathbf{N}B_{\alpha \beta }.
\label{,aSb = NBab}
\end{equation}%
where the system $B_{\alpha \beta }$, which consists of the coefficients of
proportionality between $\nabla _{\alpha }\mathbf{S}_{\beta }$ and the unit
normal $\mathbf{N}$, is known as the \textit{curvature tensor}. Thanks to
equation (\ref{Sa = ,aR}), the vector quantity $\mathbf{N}B_{\alpha \beta }$%
, which we will refer to as the \textit{vector curvature tensor}, is given
in terms of the position vector $\mathbf{R}$ by the identity%
\begin{equation}
\mathbf{N}B_{\alpha \beta }=\nabla _{\alpha }\nabla _{\beta }\mathbf{R.}
\label{NBab = ,abR}
\end{equation}%
From this equation, it immediately follows that the curvature tensor $%
B_{\alpha \beta }$ is symmetric, i.e.%
\begin{equation}
B_{\alpha \beta }=B_{\beta \alpha }.
\end{equation}%
Raising the subscript $\beta $ in equation (\ref{NBab = ,abR}) and
contracting with $\alpha $ yields%
\begin{equation}
\mathbf{N}B_{\alpha }^{\alpha }=\nabla _{\alpha }\nabla ^{\alpha }\mathbf{R.}
\label{NBAa}
\end{equation}%
The vector quantity $\mathbf{N}B_{\alpha }^{\alpha }$ is known as the 
\textit{curvature normal} by analogy with the curvature normal
characteristic of a curve. In words, the above identity says that the 
\textit{curvature normal is the surface Laplacian of the position vector}.
The invariant $B_{\alpha }^{\alpha }$ is known as the \textit{mean curvature}%
.

The covariant derivative of the unit normal $\mathbf{N}$ is given by the 
\textit{Weingarten equation}%
\begin{equation}
\nabla ^{\alpha }\mathbf{N}=-\mathbf{S}^{\beta }B_{\beta }^{\alpha }.
\label{WGE}
\end{equation}

One of the most elegant identities involving the curvature tensor is the 
\textit{Gauss equations of the surface} which read%
\begin{equation}
B_{\alpha \gamma }B_{\beta \delta }-B_{\alpha \delta }B_{\beta \gamma
}=R_{\alpha \beta \gamma \delta },  \label{GE}
\end{equation}%
where $R_{\alpha \beta \gamma \delta }$ is the Riemann-Christoffel tensor.
For a two-dimensional hypersurface in a three-dimensional Euclidean space,
the Gauss equations reduce to the form%
\begin{equation}
B_{\alpha \gamma }B_{\beta \delta }-B_{\alpha \delta }B_{\beta \gamma
}=K\varepsilon _{\alpha \beta }\varepsilon _{\gamma \delta },
\end{equation}%
where $K$ is, once again, the Gaussian curvature and $\varepsilon _{\alpha
\beta }$ is the Levi-Civita symbol. Raising the subscripts $\alpha $ and $%
\beta $ yields%
\begin{equation}
B_{\gamma }^{\alpha }B_{\delta }^{\beta }-B_{\delta }^{\alpha }B_{\gamma
}^{\beta }=K\delta _{\gamma \delta }^{\alpha \beta },  \label{BB - BB = Kd}
\end{equation}%
where $\delta _{\gamma \delta }^{\alpha \beta }=\varepsilon ^{\alpha \beta
}\varepsilon _{\gamma \delta }$ is known as the \textit{second-order delta
system}.

Note that for any second-order system $A_{\beta }^{\alpha }$ in two
dimensions we have%
\begin{equation}
A_{\gamma }^{\alpha }A_{\delta }^{\beta }-A_{\delta }^{\alpha }A_{\gamma
}^{\beta }=A\delta _{\gamma \delta }^{\alpha \beta },
\end{equation}%
where $A$ is the determinant of $A_{\beta }^{\alpha }$. Therefore, the Gauss
equations are equivalent to the statement that the Gaussian curvature equals
the determinant $B$ of $B_{\beta }^{\alpha }$, i.e.%
\begin{equation}
K=B.
\end{equation}%
In fact, an interesting generalization of the Gauss-Bonnet theorem to
arbitrary dimension is the statement that the surface integral of $B$, i.e.%
\begin{equation}
\int_{S}BdS,
\end{equation}%
depends on the topology of $S$ but not its shape.

In $n$ dimensions, write the Gauss equations%
\begin{equation}
B_{\alpha \gamma }B_{\beta \delta }-B_{\alpha \delta }B_{\beta \gamma
}=R_{\alpha \beta \gamma \delta }  \tag{\TeXButton{GE}{\ref{GE}}}
\end{equation}%
with $\alpha $ and $\beta $ as superscripts, i.e.%
\begin{equation}
B_{\gamma }^{\alpha }B_{\delta }^{\beta }-B_{\delta }^{\alpha }B_{\gamma
}^{\beta }=R_{\hspace{0.02in}\cdot \hspace{0.02in}\cdot \gamma \delta
}^{\alpha \beta }\ \ \ ,
\end{equation}%
and contract $\alpha $ with $\gamma $ and $\beta $ with $\delta $, i.e.%
\begin{equation}
B_{\alpha }^{\alpha }B_{\beta }^{\beta }-B_{\beta }^{\alpha }B_{\alpha
}^{\beta }=R_{\hspace{0.02in}\cdot \hspace{0.02in}\cdot \alpha \beta
}^{\alpha \beta }.
\end{equation}%
The invariant 
\begin{equation}
R=R_{\hspace{0.02in}\cdot \hspace{0.02in}\cdot \alpha \beta }^{\alpha \beta }
\end{equation}%
is known as the \textit{scalar curvature}. In terms of $R$, we have%
\begin{equation}
B_{\alpha }^{\alpha }B_{\beta }^{\beta }-B_{\beta }^{\alpha }B_{\alpha
}^{\beta }=R.  \label{GER}
\end{equation}%
In words, the difference between the square of the trace of the curvature
tensor $B_{\beta }^{\alpha }$ and the trace of the third fundamental form $%
B_{\beta }^{\alpha }B_{\gamma }^{\beta }$ equals the scalar curvature.

Closely related to the Gauss equations, are the Codazzi equations 
\begin{equation}
\nabla _{\alpha }B_{\beta \gamma }=\nabla _{\beta }B_{\alpha \gamma },
\end{equation}%
which, in combination with the symmetry of $B_{\beta \gamma }$, imply that
the tensor $\nabla _{\alpha }B_{\beta \gamma }$ is symmetric in all of its
subscripts. Below, we will use the following immediate consequence of the
Codazzi equations:%
\begin{equation}
\nabla _{\alpha }B_{\beta }^{\alpha }=\nabla _{\beta }B_{\alpha }^{\alpha }.
\label{C2}
\end{equation}

\section{Demonstrations of the integral identities}

\subsection{The integral $\protect\int_{S}\mathbf{N}dS$}

With the help of the two flavors of the divergence theorem -- ambient and
surface -- we can prove, in two different ways, the fact that the integral
of the unit normal $\mathbf{N}$ vanishes, i.e.%
\begin{equation}
\int_{S}\mathbf{N}dS=\mathbf{0.}  \tag{\TeXButton{IN = 0}{\ref{IN = 0}}}
\end{equation}%
Indeed, since%
\begin{equation}
\mathbf{N}=N^{i}\mathbf{Z}_{i},
\end{equation}%
we have, by the divergence theorem, that%
\begin{equation}
\int_{S}\mathbf{N}dS=\int_{S}N^{i}\mathbf{Z}_{i}dS=\int_{\Omega }\nabla ^{i}%
\mathbf{Z}_{i}dZ,
\end{equation}%
where the integrand in the last integral vanishes by the metrinilic
property. Thus, indeed,%
\begin{equation}
\int_{S}\mathbf{N}dS=\mathbf{0.}  \tag{\TeXButton{IN = 0}{\ref{IN = 0}}}
\end{equation}%
as we set out to show.

One dissatisfying aspect of this proof is the fact that it engages the
ambient space even though the elements of the integral identity 
\begin{equation}
\int_{S}\mathbf{N}dS=\mathbf{0.}  \tag{\TeXButton{IN = 0}{\ref{IN = 0}}}
\end{equation}%
are defined exclusively on the surface $S$. Therefore, we will now construct
a proof that involves only quantities defined on $S$. The presented proof
applies only to a two-dimensional hypersurface, but can be easily
generalized to arbitrary dimension.

Recall that the normal $\mathbf{N}$ is given by the identity%
\begin{equation}
\mathbf{N}=\frac{1}{2}\varepsilon ^{\alpha \beta }\mathbf{S}_{\alpha }\times 
\mathbf{S}_{\beta }.  \tag{\TeXButton{N}{\ref{N}}}
\end{equation}%
Since%
\begin{equation}
\mathbf{S}_{\alpha }=\nabla _{\alpha }\mathbf{R,}
\end{equation}%
we have%
\begin{equation}
\mathbf{N}=\frac{1}{2}\varepsilon ^{\alpha \beta }\nabla _{\alpha }\mathbf{R}%
\times \mathbf{S}_{\beta }.
\end{equation}%
By the combination of the product rule and the metrinilic property of $%
\nabla _{\alpha }$ with respect to $\varepsilon ^{\alpha \beta }$, we have%
\begin{equation}
\mathbf{N}=\frac{1}{2}\nabla _{\alpha }\left( \varepsilon ^{\alpha \beta }%
\mathbf{R}\times \mathbf{S}_{\beta }\right) -\frac{1}{2}\varepsilon ^{\alpha
\beta }\mathbf{R}\times \nabla _{\alpha }\mathbf{S}_{\beta }
\end{equation}%
Since $\mathbf{\nabla }_{\alpha }\mathbf{S}_{\beta }=\mathbf{N}B_{\alpha
\beta }$, the normal $\mathbf{N}$ is given by 
\begin{equation}
\mathbf{N}=\frac{1}{2}\nabla _{\alpha }\left( \varepsilon ^{\alpha \beta }%
\mathbf{R}\times \mathbf{S}_{\beta }\right) -\frac{1}{2}\mathbf{R}\times 
\mathbf{N}\varepsilon ^{\alpha \beta }B_{\alpha \beta }.
\end{equation}%
Next, note that since $B_{\alpha \beta }$ is symmetric and $\varepsilon
^{\alpha \beta }$ is skew-symmetric, the combination $\varepsilon ^{\alpha
\beta }B_{\alpha \beta }$ vanishes and we end up with the identity%
\begin{equation}
\mathbf{N}=\frac{1}{2}\nabla _{\alpha }\left( \varepsilon ^{\alpha \beta }%
\mathbf{R}\times \mathbf{S}_{\beta }\right)
\end{equation}%
in which the normal $\mathbf{N}$ is expressed as the surface divergence of
the combination $\varepsilon ^{\alpha \beta }\mathbf{R}\times \mathbf{S}%
_{\beta }$.

Integrating both sides of the above identity over the surface $S$, i.e.%
\begin{equation}
\int_{S}\mathbf{N}dS=\frac{1}{2}\int_{S}\nabla _{\alpha }\left( \varepsilon
^{\alpha \beta }\mathbf{R}\times \mathbf{S}_{\beta }\right) dS.
\end{equation}%
By an application of the surface divergence theorem -- recognizing that a
closed surface has no boundary $L$ -- we once again arrive at%
\begin{equation}
\int_{S}\mathbf{N}dS=\mathbf{0.}  \tag{\TeXButton{IN = 0}{\ref{IN = 0}}}
\end{equation}%
as we set out to do. Furthermore, if $S$ is a patch with a boundary $L$, we
have%
\begin{equation}
\int_{S}\mathbf{N}dS=\frac{1}{2}\int_{L}n_{\alpha }\varepsilon ^{\alpha
\beta }\mathbf{R}\times \mathbf{S}_{\beta }dL
\end{equation}%
Since the combination $n_{\alpha }\varepsilon ^{\alpha \beta }$ equals the
components $T^{\beta }$ of the unit tangent vector to the boundary $L$, we
discover that%
\begin{equation}
\int_{S}\mathbf{N}dS=\frac{1}{2}\int_{L}\mathbf{R}\times \mathbf{T}dL.
\end{equation}%
The key elements in this identity are illustrated in the following figure.%
\begin{equation}
\FRAME{itbpF}{2.2442in}{2.0358in}{0in}{}{}{sornt.png}{\special{language
"Scientific Word";type "GRAPHIC";maintain-aspect-ratio TRUE;display
"USEDEF";valid_file "F";width 2.2442in;height 2.0358in;depth
0in;original-width 13.8889in;original-height 12.583in;cropleft "0";croptop
"1";cropright "1";cropbottom "0";filename '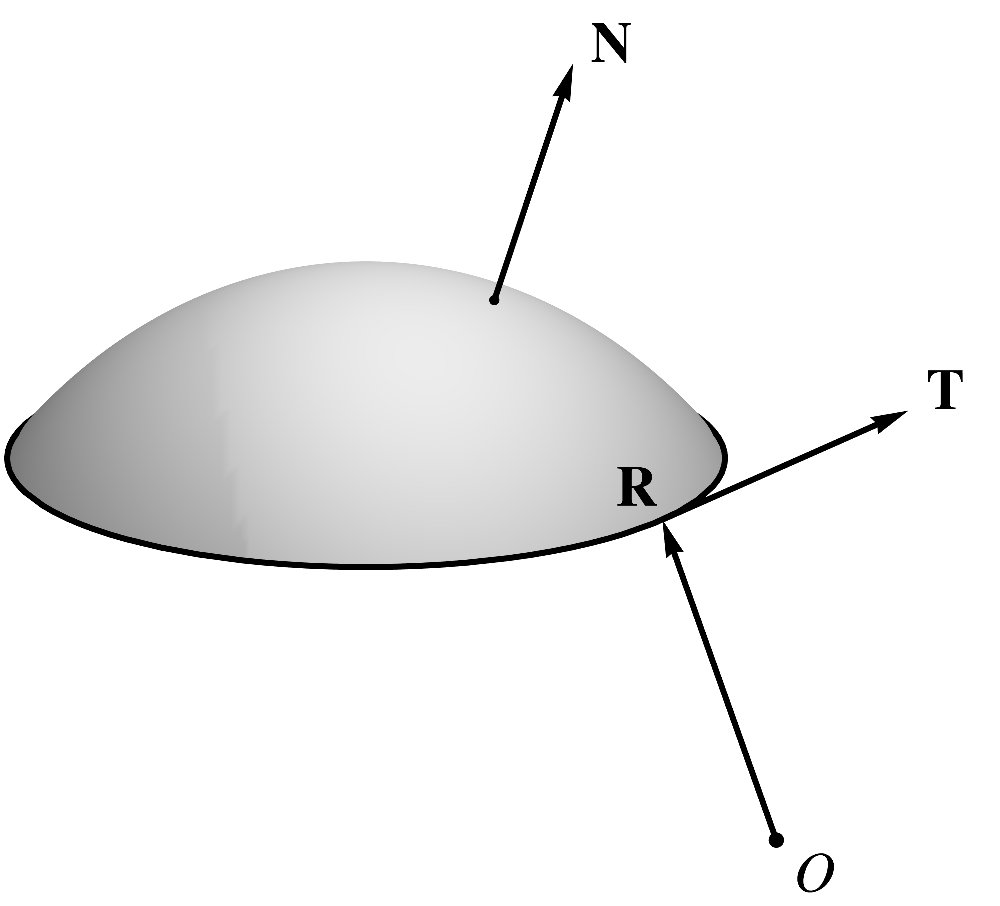';file-properties
"XNPEU";}}
\end{equation}

\subsection{The integral $\protect\int_{S}\mathbf{N}B_{\protect\alpha }^{%
\protect\alpha }dS$}

Let us now turn our attention to the identity%
\begin{equation}
\int_{S}\mathbf{N}B_{\alpha }^{\alpha }dS=\mathbf{0.} 
\tag{\TeXButton{INH
=
0}{\ref{INH = 0}}}
\end{equation}%
Recall the identity%
\begin{equation}
\mathbf{N}B_{\alpha }^{\alpha }=\nabla _{\alpha }\nabla ^{\alpha }\mathbf{R.}
\tag{\TeXButton{NBAa}{\ref{NBAa}}}
\end{equation}%
Integrating both sides over a closed surface $S$, we find%
\begin{equation}
\int_{S}\mathbf{N}B_{\alpha }^{\alpha }dS=\int_{S}\nabla _{\alpha }\nabla
^{\alpha }\mathbf{R~}dS,
\end{equation}%
where the integral on the right vanishes by the surface divergence theorem
since a closed surface has no boundary.\ In other words, we indeed find that%
\begin{equation}
\int_{S}\mathbf{N}B_{\alpha }^{\alpha }dS=\mathbf{0} 
\tag{\TeXButton{INH
=
0}{\ref{INH = 0}}}
\end{equation}%
as we set out to show.

Meanwhile, for a patch $S$ with a boundary $L$, the divergence theorem yields%
\begin{equation}
\int_{S}\mathbf{N}B_{\alpha }^{\alpha }dS=\int_{L}n_{\alpha }\nabla ^{\alpha
}\mathbf{R~}dL.
\end{equation}%
Since%
\begin{equation}
n_{\alpha }\nabla ^{\alpha }\mathbf{R}=n_{\alpha }\mathbf{S}^{\alpha }=%
\mathbf{n}
\end{equation}%
we arrive at the final identity%
\begin{equation}
\int_{S}\mathbf{N}B_{\alpha }^{\alpha }dS=\int_{L}\mathbf{n~}dL.
\end{equation}%
Two proofs of a special case of this identity can be found in \cite%
{BlackmoreSurfaceIntegral}. The contrast in complexity between those proofs
and the one presented here speaks to the effectiveness of our approach.

The key elements in the above identity are illustrated in the following
figure.%
\begin{equation}
\FRAME{itbpF}{2.386in}{2.0349in}{0in}{}{}{meancurvatureformula.png}{\special%
{language "Scientific Word";type "GRAPHIC";maintain-aspect-ratio
TRUE;display "USEDEF";valid_file "F";width 2.386in;height 2.0349in;depth
0in;original-width 13.8889in;original-height 12.583in;cropleft "0";croptop
"1";cropright "1";cropbottom "0";filename
'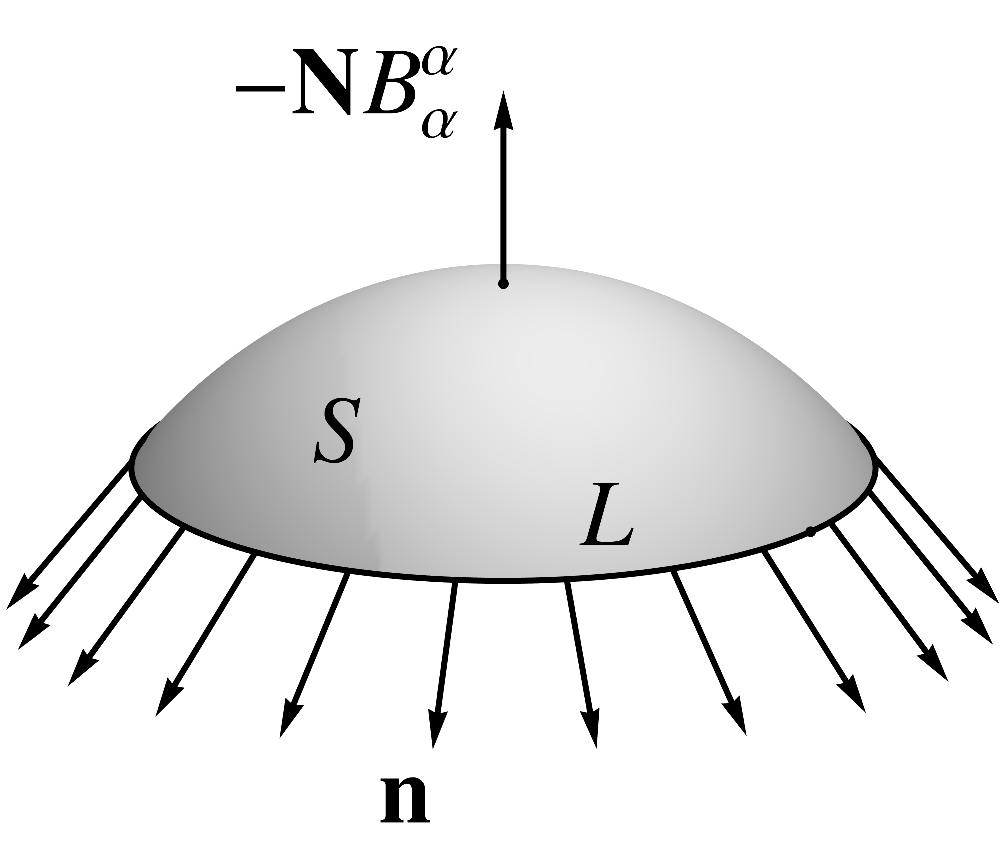';file-properties "XNPEU";}}
\end{equation}

\subsection{The integral $\protect\int_{S}\mathbf{N}KdS$}

Finally, let us demonstrate the identity%
\begin{equation}
\int_{S}\mathbf{N}KdS=\mathbf{0.}  \tag{\TeXButton{INK = 0}{\ref{INK = 0}}}
\end{equation}%
Apply the covariant derivative $\nabla _{\alpha }$ to both sides of the
Weingarten equation%
\begin{equation}
\nabla ^{\alpha }\mathbf{N}=-\mathbf{S}^{\beta }B_{\beta }^{\alpha }, 
\tag{\TeXButton{WGE}{\ref{WGE}}}
\end{equation}%
i.e.%
\begin{equation}
\nabla _{\alpha }\nabla ^{\alpha }\mathbf{N}=-\nabla _{\alpha }\left( 
\mathbf{S}^{\beta }B_{\beta }^{\alpha }\right) ,
\end{equation}%
to produce the surface Laplacian $\nabla _{\alpha }\nabla ^{\alpha }\mathbf{N%
}$ of the unit normal $\mathbf{N}$. An application of the product rule on
the right yields%
\begin{equation}
\nabla _{\alpha }\nabla ^{\alpha }\mathbf{N}=-\nabla _{\alpha }\mathbf{S}%
^{\beta }~B_{\beta }^{\alpha }-\mathbf{S}^{\beta }\nabla _{\alpha }B_{\beta
}^{\alpha }.
\end{equation}%
Since 
\begin{equation}
\nabla _{\alpha }\mathbf{S}^{\beta }=\mathbf{N}B_{\alpha }^{\beta }
\end{equation}%
and, by the Codazzi equations,%
\begin{equation}
\nabla _{\alpha }B_{\beta }^{\alpha }=\nabla _{\beta }B_{\alpha }^{\alpha },
\tag{\TeXButton{C2}{\ref{C2}}}
\end{equation}%
we have%
\begin{equation}
\nabla _{\alpha }\nabla ^{\alpha }\mathbf{N}=-\mathbf{N}B_{\alpha }^{\beta
}B_{\beta }^{\alpha }-\mathbf{S}^{\beta }\nabla _{\beta }B_{\alpha }^{\alpha
}.
\end{equation}%
Apply the "reverse" product rule to the second term on the right, i.e.%
\begin{equation}
\nabla _{\alpha }\nabla ^{\alpha }\mathbf{N}=-\mathbf{N}B_{\alpha }^{\beta
}B_{\beta }^{\alpha }-\nabla _{\beta }\left( \mathbf{S}^{\beta }B_{\alpha
}^{\alpha }\right) +\nabla _{\beta }\mathbf{S}^{\beta }B_{\alpha }^{\alpha }.
\end{equation}%
Since $\nabla _{\beta }\mathbf{S}^{\beta }=\mathbf{N}B_{\beta }^{\beta }$,
we have%
\begin{equation}
\nabla _{\alpha }\nabla ^{\alpha }\mathbf{N}=\mathbf{N}\left( B_{\alpha
}^{\alpha }B_{\beta }^{\beta }-B_{\alpha }^{\beta }B_{\beta }^{\alpha
}\right) -\nabla _{\beta }\left( \mathbf{S}^{\beta }B_{\alpha }^{\alpha
}\right) .
\end{equation}%
Now, recall that, by the Gauss equations of the surface, the quantity $%
B_{\alpha }^{\alpha }B_{\beta }^{\beta }-B_{\alpha }^{\beta }B_{\beta
}^{\alpha }$ corresponds to the scalar curvature $R$, i.e.%
\begin{equation}
B_{\alpha }^{\alpha }B_{\beta }^{\beta }-B_{\beta }^{\alpha }B_{\alpha
}^{\beta }=R.  \tag{\TeXButton{GER}{\ref{GER}}}
\end{equation}%
Thus, the surface Laplacian of $\mathbf{N}$ is given by 
\begin{equation}
\nabla _{\alpha }\nabla ^{\alpha }\mathbf{N}=\mathbf{N}R-\nabla _{\beta
}\left( \mathbf{S}^{\beta }B_{\alpha }^{\alpha }\right) .
\end{equation}%
Solving for $\mathbf{N}R$, we find%
\begin{equation}
\mathbf{N}R=\nabla _{\alpha }\left( \nabla ^{\alpha }\mathbf{N}+\mathbf{S}%
^{\alpha }B_{\beta }^{\beta }\right) .  \label{NK =}
\end{equation}%
Next, integrate both sides of the above equation over the surface $S$, i.e. 
\begin{equation}
\int_{S}\mathbf{N}RdS=\int_{S}\nabla _{\alpha }\left( \nabla ^{\alpha }%
\mathbf{N}+\mathbf{S}^{\alpha }B_{\beta }^{\beta }\right) dS.
\end{equation}%
If the surface $S$ is closed then, by the surface divergence theorem, we have%
\begin{equation}
\int_{S}\mathbf{N}RdS=\mathbf{0.}
\end{equation}%
If $S$ is a surface patch with a boundary $L$ then, by the same divergence
theorem,%
\begin{equation}
\int_{S}\mathbf{N}RdS=\int_{L}\left( n_{\alpha }\nabla ^{\alpha }\mathbf{N}%
+n_{\alpha }\mathbf{S}^{\alpha }B_{\beta }^{\beta }\right) dL\mathbf{,}
\end{equation}%
or, equivalently,%
\begin{equation}
\int_{S}\mathbf{N}RdS=\int_{L}\left( n_{\alpha }\nabla ^{\alpha }\mathbf{N}+%
\mathbf{n}B_{\alpha }^{\alpha }\right) dL\mathbf{,}
\end{equation}%
With the help of Weingarten's equation, this identity can also be written in
the form%
\begin{equation}
\int_{S}\mathbf{N}RdS=\int_{L}\left( n_{\beta }B_{\alpha }^{\alpha
}-n_{\alpha }B_{\beta }^{\alpha }\right) \mathbf{S}^{\beta }dL\mathbf{.}
\end{equation}

For a two-dimensional surface, the scalar curvature is twice the Gaussian
curvature, i.e.%
\begin{equation}
R=2K.  \label{R = 2K}
\end{equation}%
Thus, for a closed two-dimensional surface,%
\begin{equation}
\int_{S}\mathbf{N}KdS=\mathbf{0,}  \tag{\TeXButton{INK = 0}{\ref{INK = 0}}}
\end{equation}%
as we set out to show. Furthermore, for a patch with a boundary $L$, we have%
\begin{equation}
\int_{S}\mathbf{N}KdS=\frac{1}{2}\int_{L}\left( n_{\beta }B_{\alpha
}^{\alpha }-n_{\alpha }B_{\beta }^{\alpha }\right) \mathbf{S}^{\beta }dL%
\mathbf{.}
\end{equation}

\subsection{The integral $\protect\int_{S}\mathbf{R}\cdot \mathbf{N}dS$}

Let us now prove the related integral identities involving dot products with
the position vector $\mathbf{R}$. Since the following demonstrations rely on
the very same elements that we used extensively above, we will present the
demonstrations in a compressed format.

For the integral%
\begin{equation}
\int_{S}\mathbf{R}\cdot \mathbf{N}dS
\end{equation}%
we have%
\begin{eqnarray}
\int_{S}\mathbf{R}\cdot \mathbf{N}dS &=&\int_{S}\mathbf{R}\cdot N^{i}\mathbf{%
Z}_{i}dS \\
\text{(\ref{DT})} &=&\int_{\Omega }\nabla ^{i}\left( \mathbf{R}\cdot \mathbf{%
Z}_{i}\right) d\Omega \\
\text{(\ref{MN2})} &=&\int_{\Omega }\nabla ^{i}\mathbf{R}\cdot \mathbf{Z}%
_{i}~d\Omega \\
\text{(\ref{Zi = ,iR})} &=&\int_{\Omega }\mathbf{Z}^{i}\cdot \mathbf{Z}%
_{i}~d\Omega \\
\text{(\ref{ZI.Zj})} &=&\int_{\Omega }\delta _{i}^{i}~d\Omega \\
\text{(}\delta _{i}^{i}=n\text{)} &=&n\int_{\Omega }d\Omega \\
&=&nV.
\end{eqnarray}%
In the above chain, $n$ is the dimension of the ambient space and $V$ is the
volume of the enclosed domain. In summary,%
\begin{equation}
\int_{S}\mathbf{R}\cdot \mathbf{N}dS=nV,  \tag{\TeXButton{IR.N}{\ref{IR.N}}}
\end{equation}%
as we set out to show.

\subsection{The integral $\protect\int_{S}\mathbf{R}\cdot \mathbf{N}B_{%
\protect\alpha }^{\protect\alpha }dS$}

Let us now turn our attention to the integral 
\begin{equation}
\int_{S}\mathbf{R}\cdot \mathbf{N}B_{\alpha }^{\alpha }dS.
\end{equation}%
Since%
\begin{equation}
\mathbf{N}B_{\alpha }^{\alpha }=\nabla _{\alpha }\nabla ^{\alpha }\mathbf{R,}
\tag{\TeXButton{NBAa}{\ref{NBAa}}}
\end{equation}%
we have, by the product rule,%
\begin{equation}
\mathbf{R}\cdot \mathbf{N}B_{\alpha }^{\alpha }=\nabla _{\alpha }\left( 
\mathbf{R}\cdot \nabla ^{\alpha }\mathbf{R}\right) -\nabla _{\alpha }\mathbf{%
R}\cdot \nabla ^{\alpha }\mathbf{R.}
\end{equation}%
Thus, by the surface divergence theorem,%
\begin{equation}
\int_{S}\mathbf{R}\cdot \mathbf{N}B_{\alpha }^{\alpha }dS=-\int_{S}\nabla
_{\alpha }\mathbf{R}\cdot \nabla ^{\alpha }\mathbf{R}dS.
\end{equation}%
Continuing,%
\begin{eqnarray}
\int_{S}\mathbf{R}\cdot \mathbf{N}B_{\alpha }^{\alpha }dS &=&-\int_{S}\nabla
_{\alpha }\mathbf{R}\cdot \nabla ^{\alpha }\mathbf{R}dS \\
\text{(\ref{Sa = ,aR})} &=&-\int_{S}\mathbf{S}_{\alpha }\cdot \mathbf{S}%
^{\alpha }dS \\
\text{(\ref{SA.Sb})} &=&-\int_{S}\delta _{\alpha }^{\alpha }dS \\
\text{(}\delta _{\alpha }^{\alpha }=n-1\text{)} &=&-\left( n-1\right)
\int_{S}dS \\
&=&-\left( n-1\right) A.
\end{eqnarray}%
In the above chain of identities, $n$ is, once again, the dimension of the
ambient space and $A$ is the area of the surface patch $S$. In summary,%
\begin{equation}
\int_{S}\mathbf{R}\cdot \mathbf{N}B_{\alpha }^{\alpha }dS=-\left( n-1\right)
A,  \tag{\TeXButton{IR.NH}{\ref{IR.NH}}}
\end{equation}%
as we set out to show.

\subsection{The integral $\protect\int_{S}\mathbf{R}\cdot \mathbf{N}KdS$}

For greater generality, let us perform our analysis for a hypersurface in $n$
dimensions where the scalar curvature $R$ takes the place of $K$ -- or, more
precisely, of $2K$.

Recall that the expression for $\mathbf{N}R$ in divergence form reads 
\begin{equation}
\mathbf{N}R=\nabla _{\alpha }\left( \nabla ^{\alpha }\mathbf{N}+\mathbf{S}%
^{\alpha }B_{\beta }^{\beta }\right) .  \tag{\TeXButton{NK
=}{\ref{NK =}}}
\end{equation}%
Thus, by the "reverse" product rule,%
\begin{equation}
\mathbf{R}\cdot \mathbf{N}R=\nabla _{\alpha }\left( \mathbf{R}\cdot \left(
\nabla ^{\alpha }\mathbf{N}+\mathbf{S}^{\alpha }B_{\beta }^{\beta }\right)
\right) -\nabla _{\alpha }\mathbf{R}\cdot \left( \nabla ^{\alpha }\mathbf{N}+%
\mathbf{S}^{\alpha }B_{\beta }^{\beta }\right)
\end{equation}%
and, by the surface divergence theorem,%
\begin{equation}
\int_{S}\mathbf{R}\cdot \mathbf{N}RdS=-\int_{S}\nabla _{\alpha }\mathbf{R}%
\cdot \left( \nabla ^{\alpha }\mathbf{N}+\mathbf{S}^{\alpha }B_{\beta
}^{\beta }\right) dS
\end{equation}%
Continuing,%
\begin{eqnarray}
\int_{S}\mathbf{R}\cdot \mathbf{N}RdS &=&-\int_{S}\nabla _{\alpha }\mathbf{R}%
\cdot \left( \nabla ^{\alpha }\mathbf{N}+\mathbf{S}^{\alpha }B_{\beta
}^{\beta }\right) dS \\
\text{(\ref{Sa})} &=&-\int_{S}\mathbf{S}_{\alpha }\cdot \left( \nabla
^{\alpha }\mathbf{N}+\mathbf{S}^{\alpha }B_{\beta }^{\beta }\right) dS \\
\text{(\ref{WGE}),(\ref{SA.Sb})} &=&-\int_{S}\left( -\delta _{\alpha
}^{\beta }B_{\beta }^{\alpha }+\delta _{\alpha }^{\alpha }B_{\beta }^{\beta
}\right) dS \\
\text{(}\delta _{\alpha }^{\beta }B_{\beta }^{\alpha }=B_{\beta }^{\beta }%
\text{)} &=&-\int_{S}\left( \delta _{\alpha }^{\alpha }-1\right) B_{\beta
}^{\beta }dS \\
\text{(}\delta _{\alpha }^{\alpha }=n-1\text{)} &=&-\left( n-2\right)
\int_{S}B_{\alpha }^{\alpha }dS.
\end{eqnarray}%
In summary,%
\begin{equation}
\int_{S}\mathbf{R}\cdot \mathbf{N}RdS=-\left( n-2\right) \int_{S}B_{\alpha
}^{\alpha }dS,
\end{equation}

For the special case of a two-dimensional hypersurface, i.e. $n=3$, we have $%
n-2=1$ and%
\begin{equation}
R=2K,  \tag{\TeXButton{R = 2K}{\ref{R = 2K}}}
\end{equation}%
and therefore%
\begin{equation}
\int_{S}\mathbf{R}\cdot \mathbf{N}KdS=-\frac{1}{2}\int_{S}B_{\alpha
}^{\alpha }dS,  \tag{\TeXButton{IR.NH}{\ref{IR.NH}}}
\end{equation}%
as we set out to show.

\bibliographystyle{abbrv}
\bibliography{Classics,PGrinfeld,Tensors}

\end{document}